\documentclass[11pt]{article}

\usepackage{amsfonts}
\usepackage{amsmath}
\usepackage{amssymb}
\usepackage{amsthm}

\newtheorem{thm}{Theorem}[section] 
\newtheorem{prop}[thm]{Proposition} 
\newtheorem{lem}[thm]{Lemma} 
\newtheorem{cor}[thm]{Corollary} 
  
\theoremstyle{plain}

\newtheorem{conj}[thm]{Conjecture}

\newtheorem{opques}[thm]{Open Problem}

\newcommand{\tilt}{\tilde{t}}

\begin{document}

\title{Self-Organized Forest-Fires near the Critical Time}

\author{J. van den Berg and R. Brouwer \\
{\small CWI, Amsterdam } \\
{\footnotesize email: J.van.den.Berg@cwi.nl; Rachel.Brouwer@cwi.nl}
}
\date{}

\maketitle

\begin{abstract}
We consider a forest-fire model which, somewhat informally, is described
as follows:
Each site (vertex) of the square lattice is either vacant or occupied by a tree.
Vacant sites become occupied at rate $1$. Further, each site is hit by lightning
at rate $\lambda$. This lightning instantaneously destroys (makes vacant) the occupied
cluster of the site.

This model is closely related to the Drossel-Schwabl forest-fire model, which has 
received much attention in the physics literature. 
The most interesting behaviour seems to occur when the lightning rate goes to zero.
In the physics literature it is believed that then the system has so-called self-organized
critical behaviour. 

We let the system start with all sites vacant and
study, for positive but small $\lambda$,
the behaviour near the `critical time' $t_c$, defined by the relation
$1- \exp(- t_c) = p_c$, the critical probability for site percolation.

Intuitively one might expect that if, for fixed $t > t_c$, 
we let simultaneously $\lambda$ tend to $0$ and $m$ to $\infty$,
the probability that some tree at distance smaller than $m$ from
$O$ is burnt before time $t$ goes to $1$. However, we show that under a percolation-like
assumption (which we can not prove but believe to be true) this intuition is false.
We compare with the case where the square lattice is replaced by the
directed binary tree, and pose some natural open problems.

\end{abstract}

\begin{section}{Introduction}
\begin{subsection}{Background and motivation}

Consider the following, informally described, forest-fire model.
(A precise description follows later in this section).
Each site of the lattice ${\mathbb Z}^d$ is either vacant or occupied by a tree.
Vacant sites become occupied at rate $1$, independently of anything else.
Further, sites are hit by lightning at rate $\lambda$, the parameter of the model.
When a site is hit by lightning, its entire occupied cluster instantaneously burns
down (that is, becomes vacant).

This is a continuous-time version of the Drossel-Schwabl model which has received much
attention in the physics literature.  See e.g. ~\cite{BJ},
\!\! ~\cite{DS}, \! ~\cite{grass}, ~\cite{SDS}, and sections in
the book by Jensen (~\cite{Jensen}). For comparison with real forest-fires see ~\cite{MMT}.
The most interesting questions are related to the asymptotic behaviour when the lightning
rate tends to $0$. It is believed that this behaviour resembles that of `ordinary' statistical
mechanics systems at criticality. In particular, it is believed
that, asymptotically, the cluster size distribution has a power-law behaviour.
Heuristic results confirming such behaviour have been given in the literature, but
the validity of some of these results is debatable (see ~\cite{grass}) and almost nothing
is known rigorously (except for the one-dimensional case).

Our goal is more modest, and we address some basic problems which, surprisingly,
have so far been practically ignored, although their solution is crucial for a beginning
of rigorous understanding of these models. We restrict to the 2-dimensional case. That is,
the forest is represented by the square lattice.
It seems to be taken for granted in the 
literature that, informally speaking, as we let
$\lambda$ tend to $0$, the 
steady-state probability
that a given site, say the origin  $O$, is vacant stays away from $0$.
But is this really obvious? (Even, is it true?). The intuitive reasoning seems to be roughly
as follows: \\ 
``If the limit of the probability to be occupied would be $1$, then the system would 
have an `infinite occupied cluster'. But that cluster would be immediately destroyed,
bringing the occupation density away from $1$:
contradiction". \\
Of course this reasoning is, mildly speaking, quite shaky and we believe that
a rigorous solution of this problem is necessary for a clear understanding of the forest-fire model.

The problems investigated in this paper are, although not the same as the one just described,
of the same spirit.
Instead of looking at the steady-state distribution, we start with all sites
vacant and look at the time $t_c$ at which,
in the modified model where there is only growth but no ignition, an infinite cluster starts
to form. Intuitive reasoning similar to that above makes plausible that, informally
speaking, for every $t > t_c$,
the probability that $O$ burns before time $t$ stays away from $0$ as $\lambda$ tends to $0$.
Continuing such intuitive reasoning then leads to the `conclusion' that, again
informally speaking, if we take $m$ sufficiently large and
replace the above event by the event $\{$Some vertex at distance $\leq m$ from $O$ burns before
time $t \}$, the corresponding probability will be, as $\lambda$ tends to $0$, as close
to $1$ as we want. We relate this to problems which are closer to ordinary percolation.
In particular we show that,
under a percolation-like assumption (which we believe to be true), the above `conclusion'
is false.
We hope our results will lead to further research and clarification of the above problems.

\end{subsection}
\begin{subsection}{Formal statement of the problems}
So far, we have not defined our model precisely yet.
We now give this more precise description, formulate 
some of the above mentioned problems more formally, and introduce
much of the terminology used in the rest of this paper.

We work on the square lattice,
i.e. the graph of which the set of sites (vertices) is  $\mathbb{Z}^2$, and
where two vertices $(i,j)$ and $(k,l)$ share an edge if $|i-k| + |j-l| = 1$.
To each site we assign two Poisson clocks: one (which we call the `growth
clock') having rate $1$, and the other (the `ignition clock') having rate
$\lambda$.
All Poisson clocks behave independently of each other.
A site can be occupied by a tree or vacant. These states are denoted by $1$ and
$0$ respectively. Initially all sites are vacant. We restrict ourselves
to a finite box $B(n) := [-n,n]^2$. (In our theorems we consider the 
behaviour as $n \rightarrow \infty$). The dynamics is as follows:
when the growth clock of a site $v$ rings,
that site becomes occupied (unless it already was occupied, in which case
the clock is ignored); when the ignition clock of a site $v$ rings, each site that has
an occupied path in $B(n)$ to $v$, 
becomes vacant instantaneously. (Note that this means that if $v$ was already vacant,
nothing happens). 
Now let $\eta_v^n(t) = \eta_v(t) \in \{0,1\}$ denote the state of site $v$ at time $t$, and
define $\eta(t) = \eta^n(t) := (\eta_v^n(t), v \in B(n))$.
Note that, for each $n$,  $(\eta^n(t), t \geq 0)$ is a finite-state (continuous-time) irreducible
Markov chain with state space $\{0,1\}^{B(n)}$. The assignment of Poisson clocks to
{\it every} site of the square lattice provides a natural coupling of the processes
$\eta^n(\cdot),  n \geq 1$ with each other, and with other processes (see below).

For $m \leq n$, we often use the informal phrase
 ``$\eta^n$ has a fire in $B(m)$ before time $t$''
for the event $\{ \exists v \in B(m) \text{ and }\exists s \leq t \text{ such that }
\eta^n_v(s^-) = 1 \text{ and } \eta^n_v(s) = 0 \}$.
Similarly, we use  ``$\eta^n$ has at least two fires in $B(m)$ before time $t$''
for the event 
$\{ \exists v, w \in B(m) \text{ and } \exists s < u  \leq t \text{ s.t. } \eta^n_v(s^-) = 
\eta^n_w(u^-) = 1 \text{ and } \eta^n_v(s) =  \eta^n_w(u) = 0 \}$. Note that we allow $v$
and $w$ to be equal.

Let $\mathcal{P}_{\lambda}$ be the measure
that governs all the underlying Poisson processes mentioned above (and hence, for all
$n$ simultaneously, the processes $\eta^n(\cdot)$).
Often, when there is no need to explicitly indicate the dependence on $\lambda$, 
or when we consider events involving the growth clocks only,
we will omit this subscript.

It is trivial that for all times $t$ and all $n, m$ the probability that
$\eta^n$ has a fire in $B(m)$ before time $t$ goes to $0$ as $\lambda \downarrow 0$,
and hence
\begin{equation}
\lim_{n \rightarrow \infty}\ \lim_{\lambda \downarrow 0}\
\mathcal{P}_{\lambda}(\eta^n \text{ has a fire in } B(m) \text{ before time } t ) = 0.
\nonumber
\end{equation}

A much more natural (and difficult!)
question is what happens when we reverse the order of the limits.
For the investigation of such questions it turns out to be very useful to consider
the modified process $\sigma(t)$ on the infinite lattice, which we obtain, loosely
speaking, if we obey the above mentioned growth clocks but ignore the ignition clocks: 
$$\sigma_v(t) = I_{\{\text{The growth clock at } v \text{ rings in } [0,t]\}},$$
where $I$ denotes the indicator function. It is clear that, for each time $t$, the 
$\sigma_v(t), \, v \in \mathbb{Z}^2$, are Bernoulli random variables
with parameter $1 - \exp(-t)$.
So, if we define $t_c$ by the relation $p_c = 1-\exp(-t_c)$, where $p_c$ is the
critical value for ordinary site percolation on the square lattice, we see that $\sigma(t)$
has no infinite occupied cluster for $t \leq t_c$ but does have an infinite cluster for
$t > t_c$. 
To illustrate the usefulness of comparison of $\eta$ with $\sigma$ (and as introduction
to more subtle comparison arguments), we show that
\begin{equation}
\limsup_{\lambda \downarrow 0}\limsup_{n \rightarrow \infty} \mathcal{P}_{\lambda}(\eta^n
\text{ has a fire in } O \text{ before time } t) \leq \theta(1-e^{-t}),
\label{simplebound}
\end{equation}
where $\theta(.)$ is the percolation function for ordinary site percolation.
The argument is as follows: Let $\widehat{C_t}(O)$ denote the occupied cluster of $0$ in
the configuration $\sigma(t)$. It is easy to see from the process descriptions above
that in order to have, for the process $\eta^n$, a fire in $0$ before time $t$, it is 
necessary (but not sufficient) that at least one of the ignition clocks in the set 
$\widehat{C_t}(O)$ has rung before time $t$. Using the independence of the different
Poisson clocks, we have

\begin{eqnarray}
&& \hspace{-10pt}\mathcal{P}_{\lambda}(\eta^n \text{ has a fire in }
O \text{ before time } t) \nonumber \\ 
&& \hspace{15pt}\leq \sum_{k=1}^{\infty}\mathcal{P}_{\lambda}(|\widehat{C_t}(O)| = k
\text{ and } \exists v \in \widehat{C_t}(O)
\text{ that has ignition before time }t) \nonumber \\
&&\hspace{20pt} + \, \mathcal{P}(|\widehat{C_t}(O)|= \infty) \nonumber \\
&& \hspace{15pt} = \sum_{k=1}^{\infty}
\mathcal{P}(|\widehat{C_t}(O)| = k)(1-e^{-\lambda tk}) + \theta(1- e^{-t}).
\nonumber
\end{eqnarray}
Note that, in the r.h.s. above,
the first term does not depend on $n$ and,
as $\lambda \rightarrow 0$, clearly goes to $0$ (by bounded
convergence). The desired result follows.

In particular, we have for each $m$ and each $t \leq t_c$,
\begin{equation}
\lim_{\lambda \downarrow 0} \lim_{n \rightarrow \infty} \mathcal{P}_{\lambda}(\eta^n
\text{ has a fire in } B(m) \text{ before time } t ) \leq |B(m)|\theta(1-e^{-t}) = 0,
\label{nofirebefore}
\end{equation}
where $|B(m)|$ denotes the number of sites in $B(m)$.
But what happens right after $t_c$?
Intuitively one might argue as follows: \\
``If the l.h.s. of \eqref{nofirebefore} is
$0$ for some $t > t_c$, then roughly speaking, the system at time $t$ looks as in ordinary
percolation with parameter $1 - \exp(-t)$, so that an infinite occupied cluster has built
up, and this cluster intersects $B(m)$ with positive probability.
But an infinite cluster has an infinite total ignition rate and hence catches fire
immediately: contradiction. Hence for each $t > t_c$ the l.h.s. of \eqref{nofirebefore} is
strictly positive." \\
As we said before, such reasoning is very shaky. Its conclusion is correct for
the directed binary tree (see Lemma \ref{treelem}). We have some inclination to believe that
the conclusion also holds for the square lattice, but prefer to formulate this as an open
problem, rather than a conjecture:


\begin{opques}
Is, for all $t > t_c$,
\begin{equation}
\label{eqcon2}
\limsup_{\lambda \downarrow 0} \limsup_{n \rightarrow \infty}
\mathcal{P}_{\lambda}(\eta^n \text{ has a fire in } O \text{ before time }t ) > 0  \, ?
\end{equation}
\label{conj2}
\end{opques}

Believing the answer to the above problem is affirmative, it is intuitively very tempting
to go further and
`conclude' that also the answer to the following problem is affirmative:
\begin{opques}
Is it true that for all $t > t_c$ and each $\varepsilon > 0$ there exists an $m$ such
that
\begin{equation}
\label{eqquest1.2}
\limsup_{\lambda \downarrow 0} \limsup_{n \rightarrow \infty}
\mathcal{P}_{\lambda}(\eta^n \text{ has a fire in } B(m) \text{ before time }t )
> 1-\varepsilon  \, ?
\end{equation}
\label{quest1.2}
\end{opques}
The intuitive (and again shaky) reasoning here is, roughly speaking, that if 
the answer to Problem \ref{conj2} is affirmative, there will be
a positive density of sites that burn before time $t$, and hence the probability
of having such a site in B(m) will tend to $1$ as $m \rightarrow \infty$.
Our main result, Theorem \ref{onefire}, indicates that the behaviour of the 
process may be considerably different from what the above intuition suggests.
\\

At this point, one could wonder whether it is really necessary to first restrict
to finite $n$, so that we have
the annoying `extra' limit $n \rightarrow \infty$ in our theorems and problem formulations:
is, for each $\lambda >0$,  the model well-defined on the infinite lattice?
For sufficiently large $\lambda$ one can easily see that this is true.
(Using domination by suitable Bernoulli
processes one can, for such $\lambda$, make a standard graphical construction).
However, for our investigation (where we let $\lambda$ approach $0$) this is of no use.
In Section 4 we consider a slightly modified process that {\em is} defined on the infinite
lattice. In this modified process occupied clusters with size larger than or equal to $L$
(the parameter of the model) are instantaneously removed. For that model we have results very
similar to those for the original one.

In this paper we will assume knowledge of some `classical' results in 2-dimensional percolation,
in particular the standard RSW-type results (see 
\!\! ~\cite{Gr}, Chapter 11).
\label{intro}
\end{subsection}
\end{section}
\begin{section}{Statement of the main results}
\begin{subsection}{A percolation-like critical value}
In this subsection we define a percolation-like critical value, denoted by $\hat{\delta}_c$, which plays a
major role in the statement of our main results.

\smallskip\noindent
First some notational remarks. Recall that $p_c$ denotes the critical probability for 
site percolation on the square lattice. The product measure with density $p$
will be denoted by $P_p$. 
The event that there is an occupied path from a set $V$ to a set $W$ is denoted by
$\{V \leftrightarrow W\}$.
Let $n$ be a positive integer, and consider the box $B = [0,4 n] \times [0, 3 n]$. 
By the boundary of $B$, denoted by $\partial B$, we mean the set of those sites in
$B$ that have a neighbour in the complement of $B$.

We are now ready to define $\hat{\delta}_c$.
Let $\delta \in [0,1]$.
Suppose  the sites of $B$ are, independently of each other, occupied with
probability $p_c$
and vacant with probability $1-p_c$.
Next, informally, we destroy the occupied cluster of the boundary. That is, each vertex in $B$ 
that initially had an occupied path to the boundary of $B$ is made vacant. Finally, in the resulting configuration, 
each vacant site (that is, each site that initially was vacant, or that was initially occupied but made vacant 
by the above destruction step) is, independently, made occupied with probability $\delta$.
It is straightforward to see that in the final configuration a site $v \in B$ is occupied with probability
$$p_c - P_{p_c}(v \leftrightarrow \partial B) +
(1 - p_c + P_{p_c}(v \leftrightarrow \partial B)) \, \delta.$$
If we let $n$ grow and choose $v$ further
and further away from $\partial B$, this clearly converges to $p_c + (1-p_c) \delta$.
Although this is larger than
$p_c$, the final configuration has complicated spatial dependencies and therefore it is
not clear whether, in the bulk, it is
`essentially supercritical'. In particular, let $A$ be the box $[n, 3 n] \times [n, 2 n]$,  
and consider the probability $p_n(\delta)$ that 
the final configuration has an occupied vertical crossing of $A$. (As is well-known, in ordinary
supercritical percolation the probability of such event goes to $1$ as $n \rightarrow \infty$).
It is clear that $p_n(\delta)$ is increasing in $\delta$, and we define

\begin{equation}
\hat{\delta}_c = \sup \ \{ \delta : p_n(\delta) \text{ is bounded away from } 1,
\text{ uniformly in }n\}.
\label{dedef}
\end{equation}

\begin{conj} $\hat{\delta}_c > 0$.
\label{conj}
\end{conj}

In spite of serious attempts no proof (or disproof)
of this conjecture has been found yet. It is supported by simulation results but,
since the box size our simulations could handle is
limited, one has to be careful with interpreting such results.

Conjecture \ref{conj} is very similar in nature to, and `somewhat' weaker than (see the
discussion below),
Conjecture 3.2 in \!\! ~\cite{vdB&B}.
There we proved, among other results, that assumption of that conjecture yields,
informally speaking,
the non-existence of a process on the square lattice, starting with all sites vacant,
where (as in our model) vacant sites always become occupied at rate $1$,
and where infinite occupied clusters
instantaneously become vacant. Such a non-existence result, although theoretically 
interesting, looks
somewhat esoteric. In the present paper we show that the conjecture also has remarkable 
consequences for the `concrete' and natural
forest-fire models $\eta(\cdot)$.
Conjecture \ref{conj} is weaker than the above mentioned 
conjecture in \!\! ~\cite{vdB&B}, in the
sense that we can prove that the correctness of the latter implies that of the former
but we don't know how to prove the reverse implication.
Since the weaker form is sufficient for our purposes (here as well as in \!\! ~\cite{vdB&B}), we 
decided to present that form.

\end{subsection}
\begin{subsection}{The main results}
Recall the definition of $\hat{\delta}_c$ in \eqref{dedef}.
We are now ready to state our main results:

\begin{thm}
If $\hat{\delta}_c > 0$, there exists a $t> t_c$ such that for all $m$,
\begin{equation} 
\liminf_{\lambda \downarrow 0}\ \liminf_{n \rightarrow \infty}
\mathcal{P}_{\lambda}(\eta^n \text{ has a fire in } B(m) \text{ before time } t) \leq 1/2.
\end{equation}
\label{onefire}
\end{thm}

The key to Theorem \ref{onefire} is the following proposition (which is also interesting in
itself):
\begin{prop}
If $\hat{\delta}_c >0$, there exists a $t > t_c$ such that for all $m$,
\begin{equation}
\lim_{\lambda \downarrow 0}\ \limsup_{n \rightarrow \infty}\ \mathcal{P}_{\lambda}(\eta^n \text{ has at least } 
2 \text{ fires in } B(m) \text{ before time }t) = 0. 
\end{equation}
\label{no2fires}
\end{prop}
The proofs of the above proposition and theorem are given in Section 3.

\end{subsection}
\end{section}
\begin{section}{Proofs}
The proof of our main theorem (Theorem \ref{onefire}) depends heavily on Proposition \ref{no2fires}.
For the proof of the proposition we need two auxiliary models.
One of these, the `pure growth' model
$\sigma(t)$, was already introduced in Section 1. 
The other, which has the same growth mechanism but where 
removal of trees
takes place at time $t_c$ only, is described below.

\begin{subsection}{Removal at $t_c$ only}
Let $I$ denote the set of all positive even
integers $i$ 
and consider the annuli $A_i := B(5\cdot 3^i) \backslash B(3^i), i \in I$. 
Note that these annuli are pairwise disjoint. 
In the process we are going to describe,
again every site can be vacant (have value $0$) or occupied (value $1$).
By a `surrounding $i$ cluster' we will mean an occupied circuit $C$ around $0$ in
the annulus $A_i$,
together with all occupied paths in $A_i$ that contain a site in $C$.
The process is completely determined by the Poisson
growth clocks introduced in
Section 1, in the following way.
Initially each site is vacant.
Whenever the growth clock of a site rings,
the site becomes occupied. (When it already is occupied, the clock is ignored).
Destruction ($1 \rightarrow 0$ transitions) only takes place at time $t_c$:
at that time, for each positive even integer $i$,  each `surrounding $i$ cluster'
is instantaneously made vacant. After $t_c$ the growth mechanism proceeds as before.
Let $\xi_v(t)$ denote the value of site $v$ at time $t$.
Earlier in this paper we mentioned an obvious but useful relation (comparison)
between the  pure growth process $\sigma(\cdot)$ and 
the forest-fire process $\eta(\cdot)$. There is a also a
useful relation between $\xi(\cdot)$ and $\eta(\cdot)$, but its statement and
proof are less straightforward (see Lemma \ref{lemcomp} in Section 3.2).
Another lemma involving the process $\xi(\cdot)$ that will be important for us is
the following.

\begin{lem}
If $\hat{\delta}_c > 0$ there exist ${\gamma} < 1$ and $\varepsilon > 0$ such that
for all $i \in I$,
\begin{equation}
\mathcal{P}(\partial B(3^i) \rightarrow \partial B(3\cdot 3^i)
\text{ in the configuration }
\xi(t_c + \varepsilon)) < {\gamma}.
\end{equation}
\label{alreadyproven}    
\end{lem}

The proof of this lemma is very similar to that of 
Lemma 3.4 of\!\! ~\cite{vdB&B} (The $p_n$'s we defined a few lines before \eqref{dedef}
differ from the `corresponding' $a_n$'s in \!\! ~\cite{vdB&B}, but the modifications in the
proof arising from this difference are straightforward).

\end{subsection}

\begin{subsection}{Proof of Proposition \ref{no2fires}}

Fix $m$.
Since the probability in the statement of the proposition is
monotone in $m$, we may assume that $m$ is of the form $3^l$ for some even
positive integer $l$. (So each annulus $A_i, \, i \in I$, defined in the previous subsection,
is either contained in $B(m)$ or disjoint from $B(m)$).

Let $\tau = \tau(n,m)$ be the first time that $\eta^n$ has a fire in $B(m)$; 
more precisely,
\begin{equation}
\tau :=  \inf \{t : \exists v \in B(m) \text{ s.t. }\eta^n_v(t) = 0
\text{ and } \eta^n_v(t^-) = 1\}.
\nonumber
\end{equation}

Next, define, for $1 > \lambda > 0$,

\begin{eqnarray}
K(\lambda) &:=& \frac{1}{\sqrt[3]{\lambda}} \nonumber \\  
k(\lambda) &:=& \frac{1}{\sqrt[4]{\lambda}} \nonumber \\
A(k(\lambda),K(\lambda)) &:=& B(K(\lambda)) \backslash B(k(\lambda)).
\label{Kdef}
\end{eqnarray}

Further, define the following events:
\begin{eqnarray}
B_1 = B_1(\lambda) &:=& \{ \text{ no ignitions in } B(K(\lambda)) \text{ before or at time }
\tau \} \nonumber \\
B_2 = B_2(\lambda) &:=& \{\ \sigma(t_c) \text{ has a  vacant *-circuit in }
A(k(\lambda),K(\lambda)) \}, 
\nonumber
\end{eqnarray}  
where by `*-circuit' we mean a circuit (surrounding $0$) in the matching lattice (i.e.
the lattice obtained from the square lattice by adding
the two `diagonal edges' in each face of the square lattice). 

We will use the following  relation between the forest-fire process $\eta(\cdot)$ and
the auxiliary process $\xi(\cdot)$ described in the previous subsection. 
\begin{lem} Let $\lambda \in (0,1)$. 
On $B_1 \cap B_2$ we have, for
all $t > \tau$, all $v \in B(k(\lambda)) \setminus B(m)$ and all $n$, that 
\begin{equation}
\eta^{n}_v(t) \leq \xi_v(t). 
\end{equation}
\label{lemcomp}
\end{lem}
\begin{proof} (of Lemma \ref{lemcomp}).
Suppose $B_1 \cap B_2$ holds. Take $n$, $t$ and $v$ as
in the statement of the lemma. Obviously, we may assume that $k(\lambda) > m$.
To simplify notation we will, during the proof of this lemma, omit the 
superscript $n$ from $\eta$, and the argument $\lambda$ from $k$ and $K$.
Suppose $\xi_v(t) = 0$. We have to show that then
also $\eta_v(t) = 0$.
Since $\xi_v(t) = 0$,
the growth clock of $v$ does not ring in the interval $(t_c,t]$,
and we may assume that just before $t_c$
the occupied $\xi$ cluster of $v$ surrounds $B(m)$. (Otherwise the desired conclusion follows
trivially). From the definitions of the processes it then follows that at time $t_c$
the occupied $\sigma$ cluster of $v$, which we will denote by $C$, surrounds $B(m)$.
By $B_2$ we have that $C$ is in the interior
of a vacant (that is, having $\sigma(t_c) = 0$) *-circuit in $A(k,K)$.
Clearly, $\eta \equiv 0$ on this circuit during the time interval $(0, t_c]$, which
prevents fires starting 
in its exterior to reach its interior. From this, and the event $B_1$, we
conclude that $\tau > t_c$, and that $\eta(t_c)$ and $\sigma(t_c)$ agree in the interior of
this circuit. In particular, the occupied $\eta$ cluster of $v$ at time $t_c$ equals the
above mentioned set $C$.
From $B_1$ it follows that at time $\tau$ a connected set is burnt
which contains sites in $B(m)$ as well as in the complement of $B(K)$. But then
it also contains a site in $C$ (because
$C$ surrounds $B(m)$ and lies inside $B(K)$).
So the whole set $C$, and in particular $v$, burns at some time $s \in (t_c, \tau]$. Since
the growth clock of $v$ does not ring between time $t_c$ and $t$, it follows that
indeed $\eta_v(t) = 0$. This completes the proof of Lemma \ref{lemcomp}.
\end{proof}


Now we go back to the proof of the proposition.
Assume $\hat{\delta}_c > 0$.
Choose $\varepsilon$ and $\gamma$ as in Lemma \ref{alreadyproven}. 
By \eqref{nofirebefore} it is sufficient to prove that 
\begin{equation}
\lim_{\lambda \downarrow 0}\ \limsup_{n \rightarrow \infty}\
\mathcal{P}_{\lambda}(\eta^n \text{ has at least } 2 \text{ fires in } B(m)
\text{ in } (t_c, t_c + \varepsilon)) = 0.
\label{11}
\end{equation}
Define, in addition to $B_1$ and $B_2$ above, the event
$$\tilde{B}_1 = \{ \text{ no ignitions in } B(K(\lambda)) \text{ in the time interval }
(0, t_c + \varepsilon)\}.$$
We have 
\begin{eqnarray}
&& \hspace{-40pt}\mathcal{P}(\text{at least } 2 \text{ fires in } B(m) \text{ in }
(t_c, t_c + \varepsilon)) \nonumber \\ 
&& \hspace{-10pt}\leq \mathcal{P}(\{\text{at least } 2 \text{ fires in }
B(m) \text{ in } (t_c, t_c + \varepsilon)\} \cap \tilde{B}_1 \cap B_2 ) \nonumber \\ &&
+  \mathcal{P}(\tilde{B}_1^c)+ \mathcal{P}(B_2^c). 
\label{split3}
\end{eqnarray}

Now note that $\tilde{B}_1$ does not depend on $n$, and that
\begin{equation}
\mathcal{P}(\tilde{B}_1^c) \leq  \lambda\ |B(K(\lambda))|\ (t_c +	{\varepsilon}) 
\rightarrow 0, \,\, \text{ as } \lambda \downarrow 0, 
\label{B1}
\end{equation}
by the definition of $K(\lambda)$ (see \eqref{Kdef}).

Next, note that the probability of $B_2$ does not depend on $n$ either, and that
the domination of $\eta$ by $\sigma$ gives:
\begin{equation}
\mathcal{P}(B^c_2) \leq
 \mathcal{P}\{\partial B(k(\lambda)) \leftrightarrow \partial B(K(\lambda))
\text{ in } \sigma(t_c) \} \rightarrow 0, \text{ as } \lambda \downarrow 0,
\label{B2}
\end{equation}
by a well-known result from ordinary percolation and the fact that
$K(\lambda) / k(\lambda) \rightarrow \infty$ as $\lambda \downarrow 0$.

Finally, we handle the event in the first term on the right hand side of \eqref{split3}.
Since we will take limits as $\lambda \downarrow 0$, we may restrict to $\lambda$'s
for which $k(\lambda) > m$.
Then we have the following relation between events:
\begin{eqnarray}
&& \hspace{-20pt} \left\{\text{at least } 2 \text{ fires in } B(m) \text{ in }
(t_c, t_c + \varepsilon )\right\} \  \cap \ \tilde{B}_1 \ \cap \ B_2  \nonumber \\ 
&&  = \{\tau \in (t_c, t_c+ \varepsilon) \text{ and at least $1$ fire in }B(m)
\text{ in } (\tau, t_c + \varepsilon)\} \  \cap \ \tilde{B}_1 \ \cap \ B_2   \nonumber \\
&& \subset  \{ \partial B(m) \leftrightarrow \partial B(k(\lambda)) \text{ in } \eta^n(s)
\text{ for some } s \in (\tau,t_c+ \varepsilon)\} \ \cap B_1 \ \cap B_2  \nonumber \\
&& \subset  \{\partial B(m) \leftrightarrow \partial B(k(\lambda))
\text{ in } \xi(t_c + \varepsilon)\},
\label{prevlast} 
\end{eqnarray}
where the second inclusion follows from Lemma \ref{lemcomp} (and the monotonicity of
$\xi(t)$ for $t > t_c$), and the  
first inclusion holds because, by the event $\tilde{B}_1$, fires in $B(m)$ before
time $t_c+ \varepsilon$, can only arrive from outside $B(K(\lambda))$. 
To handle the probability of the last event in \eqref{prevlast}, 
first observe that, for each $i$, the random variables
$\xi_v(t)$, $t \geq 0$, $v \in A(i)$ are completely determined by 
Poisson clocks assigned to the sites inside the annulus $A(i)$.
We use the notation $I(\lambda)$ for the set of all positive even integers $j$ with
$m < 3^j  < 5\cdot 3^j \leq k(\lambda) $.  
Since the annuli $A(i), \, i \in I$ are disjoint, we get from Lemma \ref{alreadyproven}
that
\begin{equation}
\mathcal{P}\{ \partial B(m) \leftrightarrow \partial B(k(\lambda))
\text{ in }\xi(t_c + \varepsilon)\} \leq \gamma^{|I(\lambda)|}. 
\label{last}
\end{equation}
Combining \eqref{prevlast} and \eqref{last}, and using that $k(\lambda)$, and hence
also $|I(\lambda)|$ goes to $\infty$ as $\lambda \downarrow 0$, we get 
\begin{equation}
\lim_{\lambda \downarrow 0} \limsup_{n \rightarrow \infty}
\mathcal{P}(\{\text{at least } 2 \text{ fires in } B(m) \text{ before time }
(t_c + \varepsilon)\} \  \cap \tilde{B}_1 \cap B_2) = 0.
\label{allerlaatste}
\end{equation}

Combining \eqref{split3}, \eqref{B1}, \eqref{B2} and \eqref{allerlaatste} completes
the proof of Proposition \ref{no2fires}


\end{subsection}
\begin{subsection}{Proof of Theorem \ref{onefire}}
\begin{proof}
Suppose $\hat{\delta}_c > 0$ and that
for all $t > t_c$ there exists an $m = m(t)$ such that 
\begin{equation}
\liminf_{\lambda \downarrow 0}\ \liminf_{n \rightarrow \infty}
\mathcal{P}_{\lambda}( \text{ fire in } B(m) \text{ before time } t) > 1/2.
\label{assump}
\end{equation}
We will show that this leads to a contradiction.
Choose $t$ as in Proposition \ref{no2fires}.
Now take $u \in (t_c,t)$. By \eqref{assump} there 
exist $m_0$ and $\alpha(u) > 0$ such that 
\begin{equation}
\liminf_{\lambda \downarrow 0}\ \liminf_{n \rightarrow \infty}
\mathcal{P}_{\lambda}(\text{ fire in } B(m_0) \text{ before time } u) > 1/2 + \alpha(u).
\label{lowbound}
\end{equation}
By \eqref{simplebound} (and the continuity of $\theta$) we 
can choose an $s > t_c$ with
\begin{equation}
\liminf_{\lambda \downarrow 0}\ \liminf_{n \rightarrow \infty}
\mathcal{P}_{\lambda}(\text{ fire in } B(m_0) \text{ before time } s ) \leq \alpha(u)/2.
\label{l1}
\end{equation}
By \eqref{assump} there exists an $m_1 > m_0$ such that 
\begin{equation}
\liminf_{\lambda \downarrow 0}\ \liminf_{n \rightarrow \infty}
\mathcal{P}_{\lambda}(\text{ fire in } B(m_1) \text{ before time } s ) > 1/2.
\label{l2}
\end{equation}

Clearly,
\begin{eqnarray}   
&& \mathcal{P}(\text{ fire in } B(m_0) \text{ before time } u) \nonumber \\
&& \ \ \ \leq \mathcal{P}(\text{ fire in } B(m_1) \text{ before time } s \text{ and }
\text{ fire in } B(m_0) \text{ between times } s \text{ and } u) \nonumber \\
&& \ \ \ + \mathcal{P}(\text{ fire in } B(m_0) \text{ before time } s) \nonumber \\
&& \ \ \ + \mathcal{P}(no \text{ fire in } B(m_1) \text{ before time } s).
\label{clear}
\end{eqnarray}
Now for each term in \eqref{clear} we take $\liminf_{\lambda \downarrow 0}\ \liminf_{n \rightarrow \infty}$.
Then, by Proposition \ref{no2fires} the first term on the r.h.s. will vanish. Using this, and applying
\eqref{l1} and \eqref{l2} to the second and the third term respectively, yields
\begin{equation}
\liminf_{\lambda \downarrow 0}\ \liminf_{n \rightarrow \infty}\mathcal{P}( \text{ fire in } B(m_0)
\text{ before time }u) \leq 1/2 + \alpha(u)/2, 
\nonumber
\end{equation}
which contradicts \eqref{lowbound}. This completes the proof of Theorem \ref{onefire}.
\end{proof}

\end{subsection}
\end{section}
\begin{section}{Discussion and modified models}

In the model above it was the square lattice which played the role of space. Completely analogous results can be
proved, in the same way, for 
the triangular or the honeycomb lattice.

In the following subsections we discuss some less obvious modifications of the model
(different ignition mechanism; binary tree 
instead of square lattice).

\begin{subsection}{Ignition of sufficiently large clusters}
Again we work on the square lattice.
In this model the growth mechanism is the same as before (that is, 
vacant sites become occupied at rate $1$),
but the ignition mechanism is different:  Instead of the ignition rate $\lambda$ we have 
an (integer) parameter $L$.
The ignition rule now is that whenever a cluster of size $\geq L$ occurs, it is 
instantaneously ignited and burnt
down (that is, each of its sites becomes vacant). A very pleasant feature of this model is that, since the
interactions now have finite range, it can be defined on the infinite 
lattice using a standard graphical construction.
This frees us from the necessity to first work on $B(n)$ and later
take limits as
$n \rightarrow \infty$, and thus from the annoying double limits we had in our main results. 
As before, we start at time $0$ with all sites vacant.
Let $\eta_v^{[L]}(t)$ denote the value ($0$ or $1$) of site $v$ at time $t$.
The analog of Open Problem \ref{conj2} is:
\begin{opques}
Is, for all $t > t_c$  
\begin{equation}
\label{eqcon3}
\limsup_{L \rightarrow \infty}
\mathcal{P}(\eta^{[L]} \text{ has a fire in } O \text{ before time }t ) > 0  \, ?
\end{equation}
\label{conj3}
\end{opques}
Similarly, there is a straightforward analog of Open Problem \ref{quest1.2}.
Although this modified model is seemingly simpler than the original one (so that the formulation
of the problems doesn't involve the extra limit over $n$), we think the problems are,
essentially, as hard as before.

We have, with $\hat{\delta}_c$ as before (see \eqref{dedef}), analogs
of Theorem \ref{onefire} and Proposition \ref{no2fires}.
\begin{thm}
If $\hat{\delta}_c > 0$, there exists a $t> t_c$ such that for all $m$,
\begin{equation} 
\liminf_{L \rightarrow \infty} \mathcal{P}(\eta^{[L]} \text{ has a fire in } B(m)
\text{ before time } t) \leq 1/2.
\end{equation}
\label{onefireL}
\end{thm}
\begin{prop}
If $\hat{\delta}_c >0$, there exists $t > t_c$ such that for all $m$,
\begin{equation}
\lim_{L \rightarrow \infty}\ \mathcal{P}(\eta^{[L]} \text{ has at least } 
2 \text{ fires in } B(m) \text{ before time }t) = 0. 
\end{equation}
\label{no2firesL}
\end{prop}
Theorem \ref{onefireL} follows from Proposition \ref{no2firesL} in the same way as
Theorem \ref{onefire} from Proposition \ref{no2fires}.
The proof of Proposition \ref{no2firesL} is very similar to that of Proposition \ref{no2fires}
and we only indicate the main modifications:
Instead of \eqref{Kdef} we define

\begin{eqnarray}
K_L &:=& L^{1/3}, \nonumber \\
k_L  &:=& L^{1/4}.
\nonumber
\end{eqnarray}

Next, the events $B_1$, $B_2$ are replaced by the single event 
\begin{equation}
B_3 :=  \{ \sigma(t_c) \text{ has a vacant *-circuit in } A(k_L, K_L) \},
\nonumber
\end{equation}

\noindent
and Lemma \ref{lemcomp} is replaced by the following lemma, whose proof is a straightforward
modification of that of the former.
(Of course we take $m$ as before, and $\tau= \tau(L,m)$ is now defined as the first time
that $\eta^{[L]}$ has a fire in $B(m)$).

\begin{lem}
On $B_3$ we have, for
all $t > \tau$ and all $v \in B(k_L) \setminus B(m)$ that 
\begin{equation}
\eta^{[L]}_v(t) \leq \xi_v(t). 
\end{equation}
\label{lemcompL}
\end{lem}
The proof of Proposition \ref{no2firesL} now proceeds as before.
\end{subsection}

\begin{subsection}{The binary tree}
In this subsection we consider the same dynamics as for the process $\eta$ in Sections 1-3, but now 
we take the {\em directed} binary tree instead of the square lattice.

By the infinite binary directed tree, denoted by $\mathcal{T}$, we mean the tree where one vertex (called the
root) has two edges, each other vertex has three edges, and where all edges are oriented in the direction
of the root. The root will be denoted by $O$.
By the children of a site $v$ we mean the two sites {\em from} which there is an edge
{\em to} $v$. (And we say that $v$ is the parent of these sites).
By the first generation of $v$ we mean the set of children of $v$, by the second generation
the children of the children of $v$ etc.
The subgraph of $\mathcal T$ containing $O$ and its first $n$ generations will be denoted by
$\mathcal{T}(n)$.

Let us now describe the model in detail. We work on $\mathcal{T}(n)$.
Initially all sites are vacant. As in the original (Section 1) model
vacant sites become occupied at rate $1$ and occupied sites are ignited at rate $\lambda$.
When a site $v$ is ignited,
instantaneously each 
site on the occupied path from $v$ in the direction of the root is made vacant. 
The forest-fire interpretation is not very natural here.
More natural is the interpretation in terms of a nervous system:
Replace the word `site' by `node', `occupied' by `alert', vacant by
`recovering', `ignition' by `arrival of a signal from outside the system'.
Then the above description says that
whenever an alert node $v$ receives a signal (either from a child, or from outside the system),
it immediately transmits it to its
parent (except when $v=O$, in which case it `handles' the signal itself), after which it
needs an exponentially distributed recovering time to become alert again.

As before we use  $1$ to represent an occupied (`alert') and a $0$ to 
represent a vacant (`recovering') vertex.
Let $\zeta_v(t) \in \{0,1\}$ denote the state of 
vertex $v$ at time $t$.
If we want to stress dependence on $n$ we write $\zeta_v^n(t)$. 

As in Section 1, the processes $\zeta^n(\cdot)$ can be completely described in terms
of independent Poisson growth and ignition clocks, assigned to the sites of
$\mathcal T$.

Recall that site percolation on the binary tree has critical probability $1/2$, and
percolation probability function
$\theta(p) = (2 p - 1)/ p$, for $p \geq 1/2$.
Combining this with the same arguments that led to \eqref{simplebound} 
shows that, if we first let $n$ go to $\infty$ and then $\lambda$ to
$0$, the probability that the root burns before
time  $\log 2$ goes to $0$, and, moreover, that for $t > \log 2$
\begin{equation}
\label{upboundtree}
\limsup_{\lambda \downarrow 0} \limsup_{n \rightarrow \infty}
\mathcal{P}_{\lambda}(O \text{ burns before time } t) \leq  \frac{1-2e^{-t}}{1-e^{-t}}.
\end{equation}

A nice feature of the binary tree is that we can (quite simply in fact) also prove a lower
bound (compare with Open Problem \ref{conj2} for the square lattice):

\begin{lem}
\label{treelem}
For all $t > \log 2$, 
\begin{equation}
\liminf_{\lambda \downarrow 0} \limsup_{n \rightarrow \infty}
\mathcal{P}_{\lambda}(\zeta^n  \text{ has a fire in } O  \text{ before time } t)
\geq \frac{1}{2} \,\frac{1-2e^{-t}}{1-e^{-t}}.
\end{equation}
\end{lem}
Note that this lower bound is half the upper bound \eqref{upboundtree}.

\begin{proof}
Define the functions
\begin{equation}
f^{\lambda}_n(t) := \mathcal{P}_{\lambda}(\zeta^n \text{has a fire
in } O \text{ before time }t), \,\, t > 0,
\nonumber
\end{equation}
and
\begin{equation}
g^{\lambda}_n(s, t) := f^{\lambda}_n(t) -  f^{\lambda}_n(s), \,\, 0 < s < t,
\nonumber
\end{equation}
i.e. the probability that the first time that $O$ burns is between $s$ and $t$.

Fix a $t > \log 2$ and take $\tilt \in (\log 2, t)$.
Suppose that 
\begin{equation}
\liminf_{\lambda \downarrow 0}\limsup_{n \rightarrow \infty} f^{\lambda}_n(t)
< \frac{1}{2}\frac{1-2e^{-\tilt}}{1-e^{-\tilt}}.
\label{assume0}
\end{equation}
We will show that this leads to a contradiction.
By \eqref{assume0} there exists an $\alpha  > 0$ and
a sequence $(\lambda_i, \, i = 1,2, \cdots)$, which is decreasing, converges to $0$
and has, for all $i$
\begin{equation}
\limsup_{n \rightarrow \infty} f^{\lambda_i}_n(t) <
\frac{1}{2}\frac{1-2e^{-\tilt}}{1-e^{-\tilt}} - \alpha.
\label{assume}
\end{equation}
Fix $j$ large enough such that
\begin{equation}
\label{fixj}
e^{-\lambda_j \tilt}(1+2 \alpha (1-e^{-\tilt})) > 1.
\end{equation}
The reason for this choice will become clear later. 

Observe that, if $v$ and $w$ are the children of $O$,
the processes $\zeta^{n+1}_v(\cdot)$ and $\zeta^{n+1}_w(\cdot)$ are independent
copies of $\zeta^n_O(\cdot)$ (and are also independent of the Poisson clocks at $O$).
Also observe that,
to ensure that the first fire at the root occurs between times $\tilt$ and $t$, it is
sufficient that the growth clock of $O$ rings before time $\tilt$,
no ignition occurs at the root before time $\tilt$, at least one
of its children burns between times $\tilt$ and $t$ and none of its children burns before
time $\tilt$.
Hence, by these observations,
\begin{equation}
\label{recureq}
g_{n+1}^{\lambda}(\tilt,t) \geq
(1-e^{-\tilt}) \, e^{-\lambda \tilt} \, \left[ g^{\lambda}_n(\tilt,t)^2 +
2g^{\lambda}_n(\tilt,t) \, (1-f^{\lambda}_n(\tilt) \right].
\end{equation}

Now we take $\lambda$ equal to $\lambda_j$ in \eqref{recureq}, and apply \eqref{assume} 
(noting that $f^{\lambda}_n(t) \geq f^{\lambda}_n(\tilt))$. This gives that
(with the abbreviation $g_k$ for $g_k^{\lambda_j}(\tilt,t)$,
$k=1,2, \cdots$)
for all sufficiently large $n$
\begin{eqnarray}
g_{n+1} &\geq&
(1-e^{-\tilt}) \, e^{-\lambda_j\tilt} \, 2 g_n \, (1-f_n^{\lambda_j}(\tilt))
\nonumber \\
&\geq& g_n \times 
\left[ e^{-\lambda_j \tilt} (1+2 \alpha (1-e^{-\tilt}))\right].
\label{incrfn}
\end{eqnarray}
However,
the factor behind $g_n$ in the r.h.s. of \eqref{incrfn} does not depend on $n$
and is, by \eqref{fixj},
strictly larger than $1$, so that
the sequence of $g_n$'s
`explodes': a contradiction.
Hence
\begin{equation}
\label{eqconclude}
\liminf_{\lambda \downarrow 0}\limsup_{n \rightarrow \infty} f^{\lambda}_n(t)
\geq \frac{1}{2}\frac{1-2e^{-\tilt}}{1-e^{-\tilt}}.
\end{equation}
This holds for each $\tilt \in (t_c, t)$.
Letting $\tilt \uparrow t$ in \eqref{eqconclude} completes the proof of
Lemma \ref{treelem}.
\end{proof}
By a similar `independent copies' observation as used a few lines above
\eqref{recureq} (now for all sites in the $m-$th generation of the root), 
Lemma \ref{treelem} immediately gives
the following corollary (compare with Theorem \ref{onefire} and
Proposition \ref{no2fires}):

\begin{cor}
For all $t > \log 2$, all $\varepsilon > 0$ and all $k$, there exists
$m$ such that 
\begin{equation}
\liminf_{\lambda \downarrow 0}\limsup_{n \rightarrow\infty}
\mathcal{P}_{\lambda}(\zeta^n \text{ has at least } k \text{ fires in }
\mathcal{T}(m) \text{ before time } t) > 1 - \varepsilon.
\end{equation}
\end{cor}

\end{subsection}
\end{section}
{\large\bf Acknowledgments}\\
We thank Antal J\'{a}rai, Ronald Meester and Vladas Sidoravicius for
stimulating discussions.

\end{document}